\documentclass[10pt]{amsart}
\usepackage{latexsym, amsmath,amssymb}
\usepackage[active]{srcltx}

\renewcommand{\epsilon}{\varepsilon}

\newtheorem{theorem}{Theorem}
\newtheorem{lemma}[theorem]{Lemma}
\newtheorem{corr}[theorem]{Corollary}

\newtheorem{proposition}[theorem]{Proposition}
\newtheorem{deff}[theorem]{Definition}

\newcommand{\bth}{\begin{theorem}}
\newcommand{\ble}{\begin{lemma}}
\newcommand{\bcor}{\begin{corr}}

\newcommand{\bdeff}{\begin{deff}}

\newcommand{\bprop}{\begin{proposition}}
\newcommand{\ele}{\end{lemma}}
\newcommand{\ecor}{\end{corr}}
\newcommand{\edeff}{\end{deff}}

\newcommand{\eprop}{\end{proposition}}

\newcommand{\la}{\lambda}

\newcommand{\e}{\varepsilon}

\newcommand{\supp}{\text{supp }}
\renewcommand{\Pi}{\varPi}

\renewcommand{\epsilon}{\varepsilon}

\begin{document}

\title[Lower bounds on the Hausdorff measure of nodal sets]
{Lower bounds on the Hausdorff \\ measure of nodal sets II}
\thanks{The authors were supported in part by the NSF}

\author{Christopher D. Sogge}
\address{Department of Mathematics,  Johns Hopkins University,
Baltimore, MD}
\author{Steve Zelditch}
\address{Department of Mathematics, Northwestern University, Evanston, IL}

\maketitle

\begin{abstract}
We give a very short argument showing how the main identity \eqref{2} from our earlier paper \cite{SZ} immediately leads to the  best lower bound currently known \cite{CM} for the Hausdorff measure of nodal sets in dimensions $n\ge 3$.
\end{abstract}


Let $(M,g)$ be a compact smooth Riemannian manifold of dimension $n$ and let $e_\la$ be real-valued eigenfunction of the associated Laplacian, i.e.,
$$-\Delta_ge_\la(x)=\la^2e_\la(x)$$
with frequency $\la>0$.  Recent papers have been concerned with lower bounds for the $(n-1)$-dimensional Hausdorff measure, $|Z_\la|$, of the nodal set of $e_\la$, 
$$Z_\la=\{x\in M: \, e_\la(x)=0\}$$
 in dimensions $n\ge 3$.  When $n=2$ the sharp lower bound by the frequency,
$\la \lesssim |Z_\la|$, was obtained by Br\"uning in \cite{Bruning} and independently by Yau.  For all dimensions, in the analytic case, the sharp upper
and lower bounds $|Z_\la|\approx \la$ were obtained by Donnelly and Fefferman
\cite{DF1}, \cite{DF2}.

Until recently, the best known lower bound when $n\ge3$ seems to have been
$e^{-c\la}\lesssim |Z_\la|$ (see \cite{HL}).  
Using a variation \eqref{2} of an identity of Dong~\cite{Dong}, the authors showed in \cite{SZ} that
this can be improved to be $\la^{\frac74-\frac{3n}4}\lesssim |Z_\la|$. Independently 
Colding and Minicozzi \cite{CM} obtained the more favorable lower bound
\begin{equation}\label{1}
\la^{1-\frac{n-1}2}\lesssim |Z_\la|
\end{equation}
by a different method. 
Subsequently, the first author and Hezari~\cite{HS} were also able to obtain the lower bound
\eqref{1} by an argument which was in the spirit of \cite{SZ}. The purpose of this sequel to \cite{SZ} is to show
that the lower bound \eqref{1} can also  be derived by a very small modification (indeed a simplification) of the original argument of \cite{SZ}. 

 The lower bounds of \cite{SZ,HS} are  based on the identity
\begin{equation}\label{2}
\la^2 \int_M \, |e_\la|\, dV=2\int_{Z_\la} |\nabla_g e_\la|_g \, dS,
\end{equation}
from \cite{SZ} and the (sharp) lower bound for $L^1$-norms
\begin{equation}\label{3}
\la^{-\frac{n-1}4}\lesssim \int_M|e_\la|\, dV,
\end{equation}
which was also established in \cite{SZ}. Here, $dV$ is the volume element of $(M, g)$.


The lower bound \eqref{1}  is a very simple consequence of the identity
\eqref{2} and the following lemma (which was implicit in \cite{SZ}).

\begin{lemma} \label{lem} If $\la>0$ then 
\begin{equation}\label{4.0}\|\nabla_g e_\la\|_{L^\infty(M)}\lesssim \la^{1+\frac{n-1}2}\|e_\la\|_{L^1(M)}
\end{equation}
 \end{lemma}

Indeed if we use \eqref{2}  and then apply Lemma \ref{lem}, we obtain
\begin{equation} \label{EST}\begin{array}{lll} \la^2  \int_M|e_\la|\, dV
=2\int_{Z_\la}|\nabla_ge_\la|_g\, dS &\le &  2|Z_\la| \, \|\nabla_ge_\la\|_{L^\infty(M)}\\ &&\\
&\lesssim & 2|Z_\la| \; \la^{1+\frac{n-1}2}  \|e_\la\|_{L^1(M)} , \end{array}
 \end{equation}
which of course implies \eqref{1}. 

\medskip

Lemma \ref{lem} improves the upper bound on the integral given in Lemma 1 of \cite{SZ}, and its proof is almost the 
same as  the proof of \eqref{3} in Proposition 2 of \cite{SZ}:
\medskip

\begin{proof}

For  $\rho\in C^\infty_0({\mathbb R})$ we define the $\lambda$-dependent family of operators
\begin{equation}\label{5.0}
\chi_\la f = \int_{-\infty}^\infty \rho(t) e^{-it\la}e^{it\sqrt{-\Delta_g}}f\, dt=\Hat \rho(\la-\sqrt{-\Delta_g})f
=\sum_{j=0}^\infty \Hat \rho(\la-\la_j)E_jf,
\end{equation}
on $L^2(M, dV)$ with $E_jf$ denoting the projection of $f$ onto the $j$-th eigenspace of $\sqrt{-\Delta_g}$.  Here $0=\la_0<\la_1\le \la_2\le\cdots$ are its eigenvalues, and if
$\{e_j\}_{j=0}^\infty$ is the associated orthonormal basis of eigenfunctions (i.e. $\sqrt{-\Delta_g}e_j=\la_je_j$), then
$$E_jf=\Bigl(\int_M f\, \, \overline{e_j}\, dV\Bigr) \, e_j.$$
We denote the kernel of $\chi_\la$ by $K_\la (x, y)$, i.e.
$$\chi_\la f(x) = \int_M K_\la (x, y) f(y) dV(y), \;\;\; (f \in C(M)). $$
If  the Fourier transform of  $\rho$  satisfies $\Hat \rho(0)=1$, then
$\chi_\la e_\la = e_\la$, or equivalently
$$\int_M K_\la (x, y) e_\la (y) dV(y) = e_\la (x) . $$
Thus, $K_\la$ is a reproducing kernel for $e_\la$ if $\Hat \rho(0) = 1$.

As in \S 5.1 in \cite{Soggebook}, we choose $\rho$ so that the reproducing kernel  $K_\la (x, y)$  is uniformly bounded
by $\lambda^{\frac{n-1}{2}} $ on the
diagonal as $\la \to +\infty$.  This is essential for the proof of \eqref{4.0}.
If we assume that $\rho(t)=0$ for $|t|\notin [\varepsilon/2,\varepsilon]$, with $\varepsilon>0$ being a fixed number which is smaller than the
injectivity radius of $(M,g)$, then it is proved in  Lemma 5.1.3 of \cite{Soggebook} that
\begin{equation} \label{Ka} K_\la (x, y) = \lambda^{\frac{n-1}{2}} a_\la (x, y) e^{i \lambda r(x, y) }, \end{equation}
where $a_\la(x, y)$  is bounded with bounded derivatives in $(x, y)$ and where $r(x, y)$ is the Riemannian distance
between points.  This WKB formula for $K_\la(x, y)$ is known as a parametrix and may be obtained from the
H\"ormander parametrix for $ e^{i t \sqrt{-\Delta}} $ in \cite{Ho1} or from  the Hadamard parametrix for $\cos t \sqrt{-\Delta}$.
We refer to \cite{Soggebook, Hangzhou} for the background. 

 It follows from \eqref{Ka} that 
\begin{equation} \label{K} |\nabla_g K_\la(x,y)|\le C \la^{1+\frac{n-1}2}, \end{equation}
and therefore,
\begin{align*}
\sup_{x \in M} |\nabla_g \chi_\la f(x)|&= \sup_x \Bigl| \int f(y) \, \nabla_g K_\la(x,y) \, dV\Bigr|
\\
&\le \bigl\| \nabla_g K_\la(x,y) \bigl\|_{L^\infty(M \times M)}\,  \|f\|_{L^1}
\\
&\le C\la^{1+\frac{n-1}2}\|f\|_{L^1}.
\end{align*}
To complete the proof  of the Lemma, we set $f = e_\la$ and use that
$\chi_\la e_\la=e_\la$.
\end{proof}

We note that $K_\la(x, y)$ has quite a different structure from the kernels of the  spectral projection operators
$E_{[\lambda, \lambda + 1]}  = \sum_{j: \lambda_j \in [\lambda, \lambda + 1]} E_j$  and the estimate in Lemma \ref{lem} is quite different from the sup norm estimate
in Lemma 4.2.4 of \cite{Soggebook}. Indeed, in a $\lambda^{-1}$ neighborhood of the diagonal, the spectral
projections kernel $E_{[\lambda, \lambda + 1]} (x, y)$ is of size $\lambda^{n-1}$. For instance, in the case of
the standard sphere $S^n$, the kernel of the orthogonal projection $E_k$  onto the space of spherical harmonics of
degree $k \simeq \lambda$  is the constant $E_k(x,x) = \frac{\lambda^{n-1}}{Vol(S^n)}$ on the diagonal.
 We
are able to choose the test function $\rho$  above  so that the reproducing kernel  $K_\la(x, y)$ is uniformly
of size $\lambda^{\frac{n-1}{2}}$ (as in \cite{Soggebook} \S 5.1 and \cite{soggeest})  because 
we only need it to reproduce  eigenfunctions $e_\la$ of one eigenvalue
and because it does not matter  how $K_\la$ acts on   eigenfunctions of other eigenvalues. From the viewpoint
of Lagrangian distributions, the Lagrangian manifold $\Lambda_x$  associated to both $E_{[\lambda, \lambda + 1]}(x, y)$
and $K_\la (x, y)$ for fixed $x$ is the flowout $\Lambda_x = \bigcup_{t \in \supp \rho} G^t S^*_x M \subset S^* M$  of the unit-cosphere $S^*_x M$ under the geodesic flow $G^t$. The natural projection  of $\Lambda_x$ to  $M$ has  a large singularity along
$S^*_x M$ which causes the $\lambda^{n-1}$ blowup of  $E_{[\lambda, \lambda + 1]}(x, y)$ at $y = x$, but the
projection   is a covering map for the part of $\Lambda_x$ where $t \in [\epsilon, 2 \epsilon] = \mbox{supp} \rho$. The parametrix
\eqref{Ka} reflects the fact that the test function $\rho$  cuts out all of   $\Lambda_x$ except where its projection to $M$ is 
a covering map. For futher discussion of the geometry underlying Lagrangian distributions we refer to \cite{Soggebook,
Hangzhou,Z}.

Finally, we briefly  compare the proof of \eqref{1} in this note with the estimates in \cite{SZ}:

\begin{itemize}

\item Instead of Lemma \ref{1}, the estimate  $\|\nabla_g e\|_{L^\infty(M)}\lesssim \la^{1+\frac{n-1}2}\|e_\la\|_{L^2}$ was used in \cite{SZ}. The latter estimate
is a consequence of the pointwise  local Weyl law for $|\nabla e_\la(x)|^2$. 

\item In \cite{SZ} the authors proved the lower bounds \eqref{3}
by showing
that $$\|e_\la\|_{L^\infty(M)}\lesssim \la^{\frac{n-1}2}\|e_\la\|_{L^1(M)},$$
by essentially 
the same argument as in Lemma \ref{lem}.  In the proof given in this note, \eqref{3} is not used in the proof of \eqref{1} since the factor $\|e_\la\|_{L^1(M)}$
cancels out in the left and right sides. 

\end{itemize}

\end{document}